\documentclass[a4paper,10pt]{article}

\usepackage{amsfonts,amssymb,mathrsfs,amscd}
\def \R {{\mathbb {R}}}

\def\uu{\bigsqcup}
\def\eps{\varepsilon}

\usepackage[T1]{fontenc} 
\usepackage[cp1251]{inputenc}  
\usepackage[russian]{babel}

\title{\bf Slow convergence of ergodic averages \\
for actions of amenable groups }
\author{\bf Valery V. Ryzhikov}
\date{}
\begin{document}

\maketitle
\Large

For an ergodic automorphism $T$ of the probability space $(X,m)$ and a function $f\in L_1(X,m)$
theaverages $P_Nf:=\frac 1 N \sum_{i=1}^N T^if$ converge in the norm of $L_1(X,m)$ to a constant $\int fdm$.
If $f\neq 0$, then $\|P_Nf-\int fdm\|$ cannot be $o(1/N)$, since
$\|(N+1) P_Nf/N-P_{N+1}f\|=\|f\|/(N+1) $. But for suitable weighted averages  such a rate of convergence is realized (see, for example, the work of Tong and Li \cite{Izv}).
Krengel \cite{K} showed for an ergodic automorphism that the rate of convergence of uniform averages for suitable $f\in L_1(X,m)$ is arbitrarily slow. In the present communication, a similar result is established in the case of weighted averages for flows and actions of countable amenable groups.
Modifying the arguments from \cite{R}, we find a sequence of small in measure increasingly invariant sets $V_k$ such that the weighted averages for the indicator of a set lying outside the union of $V_k$ of the deviation in norm from a constant are not majorized by a predetermined sequence tending to zero. To construct such a set $V_k$, we use Rudolph's theorem on a special representation of flows (see \cite{KSF}, Chapter 11) and the Ornstein-Weiss result \cite{OW87} (II, \S 2, Theorem 5) on actions of amenable groups (an analogue of Rokhlin's lemma).

\bf Slow convergence of ergodic averages for flows. \rm \\
In the paper by Kozlov and Treschev \cite{KT}, weighted averages of the form
$P_tf(x)=\int_R h(r) f(T_{rt}x)dr$, where $\int_R h\,dr=1$, $h\geq 0$, are considered for ergodic flows. For $f\in L_1(X,m)$ we prove the convergence $\|P_tf -\int f \, dm\|\to 0$ as $t\to\infty$. We will show the absence of a universal estimate for the rate of convergence for these and more general weighted averages.

\vspace{2mm}
\bf Theorem 1. \it Let $T_t$ be an ergodic flow on the probability space $(X,m)$ that preserves the measure $m$, and let $P_n=\int_R T_t d\nu_n t$ be a sequence of operators, where $\nu_n$ are normalized Borel measures on $\R$. For any sequence $a_n \to +0$ and a set $A'$ of positive measure, for any $\eps>0$ there exists a set $A$, $m(A'\Delta A)<\eps$, such that
$$\left|\left\{n\,:\ \left\| P_n{\bf 1}_A \ - \ m(A)\right\| \ > \ a_n\right\}\right|\ =\infty.$$ \rm

\vspace{2mm}
Fix a small $\eps>0$ and a sequence ${n(k)}$ for which $\sum_ka_{n(k)}<\eps.$
The set $A$ is constructed as follows. We find $L_k$ such that $\nu_{n(k)}([-L_{n}, L_{n}])>1-\eps a_{n(k)}.$
From the special representation of the flow (see Rudolph's theorem in \cite{KSF}) it follows immediately that there exists a set $V_k$ such that $m(V_k)=2a_{n(k)}$ and $m\left(V_k\ \ \Delta \bigcap_{-L_k\leq t \leq L_k} T_tV_k\right)>(1-\eps) m(V_k).$
We set $A=A'\setminus \bigcup_k V_k$. For the operator
$Q_k=\int_{[-L_k,L_k]} T_t d\nu_n (t)$
the function $Q_{k}{\bf 1}_A $ is equal to 0 at the intersection $\bigcap_{-L_k\leq t \leq L_k} T_tV_k$. The measure of the intersection is greater than $(1-\eps) m(V_k)$ by construction. Also, due to the choice of $L_k$,
$\|P_{n(k)}{\bf 1}_A -Q_k{\bf 1}_A\|<\eps a_{n(k)} $ is satisfied. Finally, from the above, we obtain \\
$\left\|P_{n(k)}{\bf 1}_A \ - \ m(A)\right\| \ > (1-\eps) m(V_k) - \eps a_{n(k)}> a_{n(k)},$ as required.

\bf Slow convergence ofaverages for actions of countable groups. \rm Consider the free action $\{T_g\, :\, g\in G\}$ of an infinite countable group $G$ by automorphisms $T_g$ of the standard probability space $(X,m)$.
Let $F$ and $F_j$ denote finite subsets of $G$. We say that a set $V\subset X$ is $(F, c)$-invariant for $0<c <1$ if
$ m\left(\bigcap_{g\in F} T_g V\right)> c \,m(V).$
Recall that $G$ is amenable if there exists a (Follner) sequence of finite sets $F_j\subset G$ that are asymptotically invariant under group shifts (see
\cite{OW87}). Averages of the form $\sum_{g\in F_j} T_g f/|F_j|$ over such sets $F_j$ for ergodic actions converge in the norm in $L_1(X,m)$ to a constant.

\vspace{2mm}
\bf Lemma. \it Let $c\in (0,1)$, a sequence of finite
sets $F_k\subset G$, and a measure-preserving free action $\{T_g\}$ of a countable amenable group $G$ be given. If $V_k$ are $(F_k,c)$-invariant and $\sum_km(V_k) <\infty$, then for every $\eps>0$ and a set $A'$ of positive measure there exists a set $A$ such that $m(A'\Delta A)<\eps$ and
$ \left\| \sum_{g\in F_k} w_{g,k}T_g{\bf 1}_A \ - \ m(A)\right\| \ >
c\, m(V_k) m(A),$
for all non-negative $w_{g,k}$ satisfying the condition $\sum_{g\in F_k} w_{g,k}=1$.
\rm

\vspace{2mm}
For sufficiently large $N$
we set $V=\bigcup_{k>N} V_{k},$ $A:= A'\setminus V,$ so that for $U=\cup_{g\in F_k} T_g A$ we have $m(V_k\setminus U)>c m(V_k)m(A)$. But
the function $\sum_{g\in F_k} w_{g,k}T_g{\bf 1}_A$ is equal to 0 on $V_k\setminus U$, which yields the assertion of the lemma.

A consequence of the above-mentioned result of Ornstein and Weiss for free actions of amenable groups is the \it $(c,\eps)$-property \rm of these actions: \it
for any $c, \eps \in (0,1)$ and any finite $F\subset G$ there exists an $(F, c)$-invariant set of measure $\eps$. \rm
The $(c,\eps)$-property of the actions under consideration and the lemma imply the following assertion.

\vspace{2mm}
\bf Theorem 2. \it Let a sequence $a_j\to +0$, a free
ergodic action $\{T_g\}$ of an infinite countable amenable group $G$, and a sequence of finite sets $F_j\subset G$ be given. Then for $A'$, $m(A')>0$, there exists a set $A$ arbitrarily close to it in measure such that
for $w_{g,j}\geq 0$ and $\sum_{g\in F_j} w_{g,j}=1$ the following holds

\ \ \ \ \ \ \ \ $\left|\left\{j\,:\ \left\| \sum_{g\in F_j} w_{g,j}T_g{\bf 1}_A \ - \ m(A)\right\| \ > \ a_j\right\}\right|\ =\infty.$ \rm
\\
In the case of a Follner sequence $F_j$ we obtain the desired effect of slow convergence ofaverages to a constant function. \rm

\vspace{2mm}
\bf Slow convergence almost everywhere.\rm

\vspace{2mm}
\bf Theorem 3. \it Let $T$ be an ergodic automorphism of the probability space $(X,m)$ and $f\in L_1(X,m)$ be a nonzero function. For any sequence $a_n \to +0$ and a set $A'$ of positive measure, for any $\eps>0$ there exists a set $Y$, $m(Y)>1-\eps$, such that for
$\tilde f=f{\bf 1}_Y$ for almost all $x$ we have
$$\left|\left\{n\,:\ \frac 1 n \sum_{i=0}^{n-1}\tilde f(T^ix) - \int \tilde f \, dm >\ a_n\, \right\} \right|\ =\infty.$$ \rm

The proof again uses the existence of extremely high towers $V_k$ of a given (small) measure. The choice of tower $V_{k}$ is not arbitrary, but depends on the towers $V_{1},\dots, V_{k-1}$. The set $Y$ has the form
$Y=X\setminus \bigcup_k V_k$.
Let us explain how the sets $V_{k}$ are selected.
Let $\sum_k \eps_k<\eps$, $\eps_k>0$
and let the function
$$f_k= f\, {\bf 1}_{Y_k}, \ \ Y_k=X\setminus \bigcup_{i=1}^k V_i.$$
By Birkhoff's theorem, there exists $n(k)$ such that for all
$N\geq n(k)$ we have
$$m\left( x: \left|\frac 1 {N} \sum_{i=1}^{N}f_k(T^ix) -
\int f_k dm \right|< \eps_{k}\right)\ > \ 1-\eps_{k}.$$ We fix such $n(k)$ with the additional condition $\eps_{k}>a_{n(k)}$. We choose a tower
$$V_{k}=\uu_{i=1}^{h_{k}}T^iB_{k}$$
provided that $h_{k}\gg \ n(k)$ and $ m(V_{k})=3\eps_{k}.$
When choosing sufficiently large $h_{k}$ (we are free to choose arbitrarily high towers of a given measure), we have the inequality
$$m\left( \ x: \ \ \left|\frac 1 {n(k)} \sum_{i=1}^{n(k)} \tilde f(T^ix) - \int \tilde f dm \right|\ > \eps_{k} \, \right)>\ \ 1-\eps_{k}.$$
And this implies the assertion of the theorem, since $\eps_k>a_{n(k)}$.

\vspace{3mm}
The author thanks  S.V. Tikhonov, J.-P. Thouvenot and  B. Weiss for useful comments.

\vspace{20mm}
\bf  Медленная сходимость эргодических средних \\
для действий аменабельных  групп.  \rm

\vspace{5mm}
Для эргодического автоморфизма $T$ вероятностного  пространства $(X,m)$ и функции $f\in L_1(X,m)$  
 средние $P_Nf:=\frac 1 N \sum_{i=1}^N T^if$ сходятся по норме $L_1(X,m)$   к константе $\int fdm$.
 Если $f\neq 0$, то $\|P_Nf-\int fdm\|$ не может  быть $o(1/N)$, так  как 
$\|(N+1) P_Nf/N-P_{N+1}f\|=\|f\|/(N+1) $.  Но для подходящих  весовых средних такая скорость сходимости реализуется  (см., например, работу  Тонга и  Ли \cite{Izv}).  
Кренгелем \cite{K}  для  эргодического автоморфизма  показана   сколь угодно медленная скорость сходимости равномерных  средних для подходящих $f\in L_1(X,m)$.    В предлагаемом  сообщении аналогичный результат устанавливается в  случае  весовых усреднений  для  потоков  и действий счетных аменабельных групп.
Модифицируя  рассуждения  из \cite{R}, мы найдем  последовательность  малых по мере все более инвариантных  множеств $V_k$ таких, что взвешенные   средние для индикатора множества, лежащего вне  объединения $V_k$ отклонения  по норме от константы не мажорируются  заданной наперед последовательностью, стремящейся к нулю.  Для построения такого набора  $V_k$  используется  теорема Рудольфа о специальном представлении потоков (см. \cite{KSF}, глава 11)  и  результат  Орнстейна-Вейса \cite{OW87} (II, \S 2, теорема 5) о действиях аменабельных групп (аналог леммы Рохлина).    

\newpage
\bf Медленная сходимость  эргодических средних для  потоков. \rm \\
В работе Козлова и Трещева \cite{KT} для эргодических потоков рассмотрены  весовые  усреднения вида  
$P_tf(x)=\int_R h(r) f(T_{rt}x)dr$, где $\int_R h\,dr=1$, $h\geq 0$. Для $f\in L_1(X,m)$ доказана  сходимость  $\|P_tf -\int f \, dm\|\to 0$ при $t\to\infty$. Нами  будет показано отсутствие   универсальной  оценки  скорости сходимости для этих и более общих весовых усреднений. 
 
\vspace{2mm}
\bf Теорема 1. \it Пусть  $T_t$ -- эргодический поток на вероятностном пространстве $(X,m)$, сохраняющий меру $m$,  и   дана последовательность  операторов $P_n=\int_R T_t d\nu_n t$, где $\nu_n$  -- борелевские нормированные меры на $\R$. Для всякой последовательности  $a_n \to +0$  и  множества $A'$ положительной меры для любого $\eps>0$ найдется множество  $A$, $m(A'\Delta A)<\eps$, такое, что 
$\left|\left\{n\,:\ \left\|  P_n{\bf 1}_A \ - \ m(A)\right\| \ > \ a_n\right\}\right|\ =\infty.$ \rm
 
\vspace{2mm}
 Фиксируем  маленькое $\eps>0$ и последовательность ${n(k)}$, для которой $\sum_ka_{n(k)}<\eps.$ 
  Множество $A$ строится следующим образом.  Находим  $L_k$ такие, что 
$$\nu_{n(k)}([-L_{n}, L_{n}])>1-\eps a_{n(k)}.$$
Из специального представления потока (см. теорему Рудольфа в \cite{KSF}) непосредственно  вытекает существование множества $V_k$ такого, что   $m(V_k)=2a_{n(k)}$ и выполнено  $m\left(V_k\ \ \Delta \bigcap_{-L_k\leq t \leq L_k} T_tV_k\right)>(1-\eps) m(V_k).$ 
  Положим $A=A'\setminus \bigcup_k V_k$. Для оператора  
 $Q_k=\int_{[-L_k,L_k]} T_t d\nu_n (t)$
функция $Q_{k}{\bf 1}_A $ равна 0 на пересечении $\bigcap_{-L_k\leq t \leq L_k} T_tV_k$.  Мера пересечения  больше $(1-\eps) m(V_k)$ по построению. Также в силу выбора $L_k$ выполнено 
 $\|P_{n(k)}{\bf 1}_A  -Q_k{\bf 1}_A\|<\eps a_{n(k)} $. Окончательно из сказанного выше  получим \\
$\left\|P_{n(k)}{\bf 1}_A \ - \ m(A)\right\| \ > (1-\eps) m(V_k) - \eps a_{n(k)}> a_{n(k)},$ что и требовалось.

\bf Медленная сходимость средних для действий счетных групп. \rm   Рассмотрим  свободное действие $\{T_g\, :\, g\in G\}$ бесконечной счетной группы $G$  автоморфизмами   $T_g$ стандартного вероятностного пространства $(X,m)$. 
 Пусть  $F$ и $F_j$  обозначают конечные подмножества в $G$.  Говорим, что множество $V\subset X$ при $0<c <1$ является $(F, c)$-инвариантным, если 
$  m\left(\bigcap_{g\in F} T_g V\right)> c \,m(V).$
Напомним, что $G$ аменабельна, если найдется  (фёльнеровская) последовательность конечных множеств $F_j\subset G$, асимптотически инвариантных относительно групповых сдвигов (см.
\cite{OW87}). Усреднения вида $\sum_{g\in F_j} T_g f/|F_j|$   по таким множествам $F_j$ для эргодических действий сходятся по норме в $L_1(X,m)$ к константе.

\vspace{2mm} 
\bf Лемма.   \it    Пусть дано $c\in (0,1)$, последовательность  конечных
множеств   $F_k\subset G$ и    сохраняющее меру  свободное действие $\{T_g\}$ счетной аменабельной группы $G$. Если $V_k$ являются   $(F_k,c)$-инвариантными и  $\sum_km(V_k) <\infty$, то для всяких  $\eps>0$ и  множества $A'$ положительной меры найдется множество  $A$ такое, что  $m(A'\Delta A)<\eps$  и 
$ \left\| \sum_{g\in F_k} w_{g,k}T_g{\bf 1}_A \ - \ m(A)\right\| \ > 
c\, m(V_k) m(A),$ 
для всяких  неотрицательных  $w_{g,k}$,  удовлетворяющих условию $\sum_{g\in F_k} w_{g,k}=1$.  
\rm 

\vspace{2mm} 
 Для достаточно большого  $N$ 
положим  $V=\bigcup_{k>N} V_{k},$  $A:= A'\setminus V,$ чтобы для  $U=\cup_{g\in F_k} T_g A$  выполнялось   $m(V_k\setminus U)>c m(V_k)m(A)$. Но
функция $\sum_{g\in F_k} w_{g,k}T_g{\bf 1}_A$ равна 0 на $V_k\setminus U$, что дает  утверждение леммы.

Cледствием упомянутого результа Орнстейна и Вейса для свободных действий аменабельных групп     является \it $(c,\eps)$-свойство \rm этих действий: \it
  для  всяких  $c, \eps \in (0,1)$  и любого конечного $F\subset G$ найдется  $(F, c)$-инвариантное множество  меры $\eps$. \rm
Из  $(c,\eps)$-свойства рассматриваемых действий и леммы вытекает следующее утверждение.    
 
\vspace{2mm} 
\bf Теорема 2. \it Пусть заданы последовательность $a_j\to +0$, свободное 
эргодическое действие $\{T_g\}$ бесконечной счетной аменабельной группы $G$ и  последовательность конечных множеств $F_j\subset G$. Тогда для  $A'$, $m(A')>0$, найдется сколь угодно близкое к нему по мере множество $A$  такое, что 
при  $w_{g,j}\geq 0$ и $\sum_{g\in F_j} w_{g,j}=1$   выполняется 

\ \ \ \ \ \ \ \ $\left|\left\{j\,:\ \left\| \sum_{g\in F_j} w_{g,j}T_g{\bf 1}_A \ - \ m(A)\right\| \ > \ a_j\right\}\right|\ =\infty.$  
\\
В случае фёльнеровской последовательности  $F_j$  получаем нужный нам  эффект медленной сходимости средних к постоянной функции.  \rm

\vspace{2mm} 
\bf Медленная сходимость почти всюду.\rm

\vspace{2mm} 
\bf Теорема 1. \it Пусть  $T$ -- эргодический автоморфизм вероятностного пространства $(X,m)$ и   $f\in L_1(X,m)$ -- ненулевая функция.  Для всякой последовательности  $a_n \to +0$  и  множества $A'$ положительной меры для любого $\eps>0$ найдется множество  $Y$, $m(Y)>1-\eps$, такое, что  для 
$\tilde f=f\,{\bf 1}_Y$ для почти всех $x$  выполнено  
$$\left|\left\{n\,:\ \frac 1 n \sum_{i=0}^{n-1}\tilde f(T^ix) - \int \tilde f\, dm >\ a_n\, \right\} \right|\ =\infty.$$ \rm

Доказательство опять  использует наличие  чрезвычайно высоких башен $V_k$ заданной (малой) меры.  Выбор  башни $V_{k}$ не произволен, а зависит от   башен $V_{1},\dots,  V_{k-1}$.  Множество $Y$  имеет вид
$Y=X\setminus \bigcup_k  V_k$. 
Поясним, как выбираются  множества   $V_{k}$.
Пусть $\sum_k \eps_k<\eps$, $\eps_k>0$ 
 и задана  функция 
$$f_k= f\, {\bf 1}_{Y_k},  \ \ Y_k=X\setminus \bigcup_{i=1}^k  V_i.$$ 
По теореме Биркгофа найдется  $n(k)$ такое, что для всех
$N\geq n(k)$ выполнено 
  $$m\left( x:  \left|\frac 1 {N} \sum_{i=1}^{N}f_k(T^ix) - 
\int f_k dm \right|< \eps_{k}\right)\ > \ 1-\eps_{k}.$$    Фиксируем такое $n(k)$  с дополнительным условием  $\eps_{k}>a_{n(k)}$.  Выбираем   башню 
$$V_{k}=\uu_{i=1}^{h_{k}}T^iB_{k}$$
 при условии, что $h_{k}\gg \ n(k)$ и $  m(V_{k})=3\eps_{k}.$
При выборе достаточно   больших $h_{k}$ (мы же свободны в выборе сколь угодно высоких башен заданной меры) будет обеспечено неравенство  
$$m\left(\ x:  \ \left| \frac 1 {n(k)} \sum_{i=1}^{n(k)} \tilde f(T^ix) - \int \tilde f dm \right| \  > \eps_{k} \, \right)>\ \ 1-\eps_{k}.$$
А это влечет за собой утверждение теоремы, так как $\eps_k>a_{n(k)}$.

\vspace{3mm} 
Автор благодарит  Б. Вейса,  С.В. Тихонова и  Ж.-П. Тувено за  полезные замечания.

\end{document}